\documentclass[10pt]{article}

\usepackage{amsmath}
\usepackage{amsfonts}
\usepackage{amssymb}

\newtheorem{theorem}{Theorem}[section]
\newtheorem{lemma}[theorem]{Lemma}
\newtheorem{corollary}[theorem]{Corollary}
\newtheorem{remark}[theorem]{Remark}

\def\phi{{\varphi}}

\DeclareSymbolFont{AMSb}{U}{msb}{m}{n}
\DeclareMathSymbol{\N}{\mathbin}{AMSb}{"4E}
\DeclareMathSymbol{\Z}{\mathbin}{AMSb}{"5A}
\DeclareMathSymbol{\R}{\mathbin}{AMSb}{"52}
\DeclareMathSymbol{\Q}{\mathbin}{AMSb}{"51}
\DeclareMathSymbol{\I}{\mathbin}{AMSb}{"49}
\DeclareMathSymbol{\C}{\mathbin}{AMSb}{"43}

\begin{document}

\title{Optimal weighted Hardy-Rellich inequalities on $H^2\cap H^{1}_{0}$}
\author{  Amir Moradifam \thanks{This work is supported by a Killam Predoctoral
Fellowship, and is part of the author's PhD dissertation in
preparation under the supervision of N. Ghoussoub.}\\
\small Department of Mathematics,
\small University of British Columbia, \\
\small Vancouver BC Canada V6T 1Z2 \\
\small {\tt a.moradi@math.ubc.ca}
\\
%\today\\
%\date{January 20, 2005}\\
} \maketitle

 \begin{abstract}

We give necessary and sufficient conditions  on a pair of positive radial %%@
functions  $V$ and $W$ on a ball $B$ of radius $R$  in $R^{n}$, $n
\geq 1$, so that the following inequalities hold
\begin{equation*} \label{two}
\hbox{$\int_{B}V(x)|\nabla u |^{2}dx \geq  \int_{B} W(x) %%@
u^{2}dx+b\int_{\partial B}u^2 ds$\ \ for all u $\in H^1(B)$,}
\end{equation*}
and
\begin{equation*} \label{two}
\hbox{$\int_{B}V(x)|\Delta u |^{2}dx \geq  \int_{B} W(x)|\nabla  %%@
u|^{2}dx+b\int_{\partial B}|\nabla u|^2 ds$ \ \ for all u $\in
H^2(B)$.}
\end{equation*}
Then we present various classes of optimal weighted Hardy-Rellich
inequalities on $H^{2}\cap H^{1}_{0}$. The proofs are based on
decomposition into spherical harmonics. These types inequalities are
important in the study of fourth order elliptic equations with
Navier boundary condition and systems of second order elliptic
equations.
 \end{abstract}

\section{Introduction}

Let $\Omega$ be a  smooth bounded domain in $\R^n$ and $0 \in
\Omega$. Let us recall that the classical Hardy-Rellich inequality
assets that
\begin{equation}\label{HR}
\int_{\Omega}|\Delta u|^{2}dx \geq \frac{n^2(n-4)^2}{16}
\int_{\Omega}\frac{u^{2}}{|x|^{4}}dx,\ \ \hbox{for}\ \ u \in
H^2_{0}(\Omega),
\end{equation}
where the constant appearing in the above inequality is the best
constant and it is never achieved in $H^2_{0}$. Recently there has
been a flurry of activity about possible improvements of the
following type

\begin{equation}\label{gen-rellich.0}
\hbox{If $n\geq 5$\quad then \quad $\int_{\Omega}|\Delta u|^{2}dx -
\frac{n^2(n-4)^2}{16} \int_{\Omega}\frac{u^{2}}{|x|^{4}}dx\geq
 \int_{\Omega} W(x)u^{2}dx$ \quad for $u  \in  H^{2}_{0}(\Omega)$},
\end{equation}
as well as
\begin{equation}\label{sec-hardy.0}
\hbox{If $n\geq 3$\quad then \quad $\int_{\Omega}|\Delta u |^{2}dx - C(n) %%@
\int_{\Omega}\frac{|\nabla u|^{2}}{|x|^{2}}dx\geq \int_{\Omega} V(x)|\nabla u|^{2}dx$ \quad for all $u \in %%@
H^{2}_{0}(\Omega)$,}
\end{equation}
where $V, W$ are certain explicit radially symmetric potentials of order lower than %%@
$\frac{1}{r^2}$ (for $V$) and  $\frac{1}{r^4}$ (for $W$) (see
\cite{AGS}, \cite{B}, \cite{C}, \cite{D1}, \cite{DGLV}, \cite{GM2},
and \cite{TZ}.

The inequality $(\ref{HR})$ was first proved by Rellich \cite{R} for
$u \in H^{2}_{0}(\Omega)$ and then it was extended to functions in
$H^{2}(\Omega)\cap H^{1}_{0}(\Omega)$ by Donal et al. in
\cite{DGLV}. So far most of the results about improved Hardy-Rellich
inequalities and the inequalities of the form (\ref{sec-hardy.0})
are proved for $u \in H^{2}_{0}(\Omega)$ (see \cite{C}, \cite{GM2},
and \cite{TZ}). The goal of this paper is to provide a general
approach to prove optimal weighted Hardy-Rellich inequalities on
$H^{2}(\Omega)\cap H^{1}_{0}(\Omega)$ and inequalities of type
(\ref{sec-hardy.0}) on $H^{2}(\Omega)$ which are important in the
study of fourth order elliptic equations with Navier boudary
condition and systems of second order elliptic equations (see
\cite{Mo}).

\quad We start -- in section 2 -- by giving necessary and sufficient conditions  on positive %%@
radial functions  $V$ and $W$ on a ball $B$  in $R^{n}$,   so that the following inequality %%@
holds for some $c>0$ and $b<0$:
\begin{equation}\label{most.general.hardy}
\hbox{$\int_{B}V(x)|\nabla u |^{2}dx \geq c\int_{B} W(x)u^2dx+b\int_{\partial B}u^2$ for all $u \in %%@
H^{1}(B)$.}
\end{equation}
Assuming that the ball $B$ has radius $R$ and that %%@
$\int^{R}_{0}\frac{1}{r^{n-1}V(r)}dr=+\infty$,  the condition is simply that the ordinary %%@
differential equation
\begin{equation*}
\hbox{ $({\rm B}_{V,cW})$  \quad \quad \quad \quad \quad \quad \quad \quad \quad \quad \quad $y''(r)+(\frac{n-1}{r}+\frac{V_r(r)}{V(r)})y'(r)+\frac{cW(r)}{V(r)}y(r)=0$ \quad \quad \quad \quad \quad \quad \quad \quad \quad \quad \quad}
\end{equation*}
has a positive solution $\phi$ on the interval $(0, R)$ with
$V(R)\frac{\phi'(R)}{\phi(R)}=b$. As in \cite{GM2}, we shall call
such a couple $(V, W)$ a {\it Bessel pair on $(0, R)$}. The {\it
weight} of such a pair is then defined as
\begin{equation}
\hbox{$\beta (V, W; R)=\sup \big\{ c;\,  ({\rm B}_{V,cW})$ has a positive solution  on $(0, R)\big\} $.}
\end{equation}
We call $W$ a Bessel potential if $(1,W)$ is a Bessel pair. This
characterization makes an important connection between Hardy-type
inequalities and the oscillatory behavior of the above equations.
For a detailed analysis of Bessel pairs see \cite{GM2}. The above
theorem in the general form of improved Hardy-type inequalities
which recently has been of interest for many authors (see
\cite{ACR}, \cite{BV}, \cite{CKN}, \cite{CW}, \cite{Co}, \cite{D},
\cite{FT},
\cite{GGM}, \cite{VZ}, and \cite{WW}).\\

Here is the main result of this paper.

\begin{theorem}  Let $V$ and $W$ be positive radial $C^1$-functions   on $B\backslash \{0\}$, %%@
where $B$ is a ball centered at zero with radius $R$ in $\R^n$ ($n \geq 1$) such that  %%@
$\int^{R}_{0}\frac{1}{r^{n-1}V(r)}dr=+\infty$ and $\int^{R}_{0}r^{n-1}V(r)dr<+\infty$. The %%@
following statements are then equivalent:

\begin{enumerate}

\item $(V, W)$ is a Bessel pair on $(0, R)$ with $\theta:= V(R)\frac{\phi'(R)}{\phi(R)}$, where $\phi$ is the corresponding solution of
$(B_{(V,W)})$.

\item $ \int_{B}V(x)|\nabla u |^{2}dx \geq \int_{B} W(x)u^2dx+\theta \int_{\partial B}u^2 ds$ for all $u \in C^{\infty}(\bar{B})$.

\item If $\lim_{r \rightarrow 0}r^{\alpha}V(r)=0$ for some $\alpha< n-2$, then the above are %%@
equivalent to
\[
\hbox{$\int_{B}V(x)|\Delta u |^{2}dx \geq  \int_{B} W(x)|\nabla  %%@
u|^{2}dx+(n-1)\int_{B}(\frac{V(x)}{|x|^2}-\frac{V_r(|x|)}{|x|})|\nabla
u|^2dx+(\theta+(n-1)V(R))\int_{\partial B}|\nabla u|^2$,}
\]
for all radial $u \in C^{\infty}(\bar{B})$.

 \item If in addition, $W(r)-\frac{2V(r)}{r^2}+\frac{2V_r(r)}{r}-V_{rr}(r)\geq 0$ on $(0, R)$, %%@
then the above are equivalent to
\[
\hbox{$\int_{B}V(x)|\Delta u |^{2}dx \geq  \int_{B} W(x)|\nabla  %%@
u|^{2}dx+(n-1)\int_{B}(\frac{V(x)}{|x|^2}-\frac{V_r(|x|)}{|x|})|\nabla
u|^2dx+(\theta+(n-1)V(R))\int_{\partial B}|\nabla u|^2$,}
\]
for all  $u \in C^{\infty}(\bar{B})$.
\end{enumerate}
\end{theorem}

Appropriate combinations of  $4)$ and $2)$ in the above theorem and  lead  %%@
to  a myriad of Hardy-Rellich type inequalities on
$H^{2}(\Omega)\cap H^{1}_{0}(\Omega)$.

\begin{remark}
The condition
$W(r)-\frac{2V(r)}{r^2}+\frac{2V_r(r)}{r}-V_{rr}(r)\geq 0$ in the
above theorem guarantees that the minimizing sequences are radial
functions. We shall see in section 3 that even with out this
condition our approach is applicable, although  the minimizing
sequences are no longer radial functions.

\end{remark}

\begin{remark}
To see the importance and generality of the above theorem, notice
that inequalities (7) and (8) in \cite{Mo} which are the author's
main tools to prove singularity of the extremal solutions in
dimensions $n\geq 9$ (see \cite{Mo}) are an immediate consequence of
the above theorem combined with (\ref{most.general.hardy}). This
theorem will also allow us to extend most of the results about Hardy
and Hardy-Rellich type inequalities on  $C^{\infty}_{0}(\Omega)$ to
corresponding inequalities on $C^{\infty}(\bar{\Omega})$ such as
those in \cite{GM2} and \cite{TZ}.
\end{remark}

We shall show that for $-\frac{n}{2}\leq m\leq\frac{n-2}{2}$
\begin{equation}
H_{n,m}=\inf_{u \in H^2(B)\setminus\{0\}}\frac{\int_{B}\frac{|\Delta
u|^2}{|x|^2m}}{\int_{B}\frac{|\nabla u|^2}{|x|^{2m+2}}}=\inf_{u \in
H_{0}^2(B)\setminus \{0\}}\frac{\int_{B}\frac{|\Delta
u|^2}{|x|^2m}}{\int_{B}\frac{|\nabla u|^2}{|x|^{2m+2}}},
\end{equation}
and for $-\frac{n}{2}\leq m\leq\frac{n-4}{2}$
\begin{equation}
a_{n,m}=\inf_{u \in H^2(B)\cap
H^{1}_{0}(B)\setminus\{0\}}\frac{\int_{B}\frac{|\Delta
u|^2}{|x|^{2m}}}{\int_{B}\frac{u^2}{|x|^{2m+4}}}=\frac{\int_{B}\frac{|\Delta
u|^2}{|x|^{2m}}}{\int_{B}\frac{u^2}{|x|^{2m+4}}},
\end{equation}
where the constants $H_{n,m}$ and $a_{n,m}$ have been computed in
\cite{TZ} and then more generally in \cite{GM2}. For example
$a_{n,0}=\frac{n^2}{4}$ for $n\geq 5$, $a_{4,0}=3$, and
$a_{3,0}=\frac{25}{36}$.

The above general theorem also allows us to obtain improved
Hardy-Rellich inequalities on $H^2(B)\cap
H^{1}_{0}(B)$. For instance, assume $W$ is a  %%@
Bessel potential on $(0, R)$  and $\phi$ is the corresponding solution of $(B_{(1,W)})$ with $R\frac{\phi'(R)}{\phi(R)}\geq -\frac{n}{2}$. If $r\frac{W_r(r)}{W(r)}$ decreases to %%@
$-\lambda$ and $\lambda \leq n-2$, then we have for all $ H^2(B)\cap
H^{1}_{0}(B)$
  \begin{equation} \label{Rellich.1}
\int_{B}|\Delta u|^{2}dx -
\frac{n^2(n-4)^2}{16}\int_{B}\frac{u^2}{|x|^4}dx\geq %%@
\big(\frac{n^2}{4}+\frac{(n-\lambda-2)^2}{4}\big)
 \beta (W; R)\int_{B}\frac{W(x)}{|x|^2}u^2 dx.
\end{equation}

By applying (\ref{Rellich.1}) to the various examples of Bessel
functions,  we can various  improved Hardy-Rellich inequalities on $
H^2(B)\cap H^{1}_{0}(B)$. Here are some basic examples of Bessel
potentials, their corresponding solution $\phi$ of $(B_{(1,W)})$.
\begin{itemize}
\item $  W \equiv 0$ is a Bessel potential on $(0,  R)$ for any
$R>0$ and $\phi=1$.
\item $  W \equiv 1$ is a Bessel potential on $(0,  R)$ for any $R>0$, $\phi(r)=J_{0}(\frac{\mu r}{R})$, where $J_{0}$ is the Bessel function and
$z_{0}=2.4048...$ is the first zero of the Bessel function $J_0$.
Moreover $R\frac{\phi'(R)}{\phi(R)}=-\frac{n}{2}$.

\item  For  $k\geq 1$, $R>0$, let  $ W_{k, %%@
\rho} (r)=\Sigma_{j=1}^k\frac{1}{r^{2}}\big(\prod^{j}_{i=1}log^{(i)}\frac{\rho}{r}\big)^{-2}$ %%@
where the functions $log^{(i)}$ are defined  iteratively  as follows:  $log^{(1)}(.)=log(.)$ and  %%@
for $k\geq 2$,  $log^{(k)}(.)=log(log^{(k-1)}(.))$.  $ W_{k, \rho}$ is then a Bessel potential on %%@
$(0, R)$ with the corresponding solution
\[\phi_{k}=\big(\prod^{j}_{i=1}log^{(i)}\frac{\rho}{r}\big)^{-\frac{1}{2}}.\]
It is easy to see that for $\rho\geq R(
e^{e^{e^{.^{.^{e((k-1)-times)}}}}} )$ large enough we have
$R\frac{\phi_{k}'(R)}{\phi_{k}(R)}\geq-\frac{n}{2}$.

\item  For $k\geq 1$, and $R>0$,  define $\tilde W_{k; \rho} %%@
(r)=\Sigma_{j=1}^k\frac{1}{r^{2}}X^{2}_{1}(\frac{r}{R})X^{2}_{2}(\frac{r}{R})
\ldots
X^{2}_{j-1}(\frac{r}{R})X^{2}_{j}(\frac{r}{R})$ where  the functions $X_i$ are defined %%@
iteratively  as follows:
 $X_{1}(t)=(1-\log(t))^{-1}$ and for $k\geq 2$, $ X_{k}(t)=X_{1}(X_{k-1}(t))$. Then again $ \tilde %%@
W_{k, \rho}$ is a Bessel potential on $(0, R)$ with
$\phi_{k}=(X_{1}(\frac{r}{R})X_{2}(\frac{r}{R}) \ldots
X_{j-1}(\frac{r}{R})X_{k}(\frac{r}{R}))^{\frac{1}{2}}$. Moreover,
$R\frac{\phi'_{k}(R)}{\phi_{k}(R)}=-\frac{k}{2}$.

\end{itemize}

As an example, let $k\geq 1$ and choose $\rho\geq R(
e^{e^{e^{.^{.^{e(k-times)}}}}} )$ large enough so that
$R\frac{\phi'(R)}{\phi(R)}\geq -\frac{n}{2}$, where
\begin{equation}\label{log.def}
\phi=\big(
\prod^{j}_{i=1}log^{(i)}\frac{\rho}{|x|}\big)^{\frac{1}{2}}.
\end{equation}
Then we have
\begin{equation}
 \int_{B}|\Delta u(x) |^{2}dx \geq \frac{n^2(n-4)^2}{16}\int_{B}\frac{u^2}{|x|^4} %%@
dx+(1+\frac{n(n-4)}{8})\sum^{k}_{j=1}\int_{B}\frac{u^2}{|x|^4}\big( %%@
\prod^{j}_{i=1}log^{(i)}\frac{\rho}{|x|}\big)^{-2}dx,
 \end{equation}
 for all $H^2(B)\cap H^{1}_{0}(B)$ which corresponds to the result
 od Adimurthi et al. \cite{AGS}.
% and
% \begin{equation}
% \int_{B}|\Delta u(x) |^{2}dx \geq \frac{n^2(n-4)^2}{16}\int_{B}\frac{u^2}{|x|^4} %%@
%dx+(1+\frac{n(n-4)}{8})\sum^{k}_{j=1}\int_{B}\frac{u^2}{|x|^4}X^{2}_{1}(\frac{|x|}{\rho})X^{2%%@
%%@
%}_{2}(\frac{r}{\rho}) \ldots
%X^{2}_{j-1}(\frac{|x|}{\rho})X^{2}_{j}(\frac{|x|}{\rho})dx.
% \end{equation}

More generally, we show that for any  $-\frac{n}{2}\leq m<\frac{n-2}{2}$, and any $W$ Bessel potential on a ball %%@
$B_{R}\subset R^n$ of radius $R$, if for the corresponding solution
$\phi$ of ($B_{(1,W)}$) we have
$R\frac{\phi'(R)}{\phi(R)}\geq-\frac{n}{2}-m$, then the following
inequality holds for all $u \in C^{\infty}_{0}(B_{R})$
\begin{equation}\label{gm-hr.00}
\int_{B_{R}}\frac{|\Delta u|^2}{|x|^{2m}}\geq a_{n,m}\int_{B_{R}}\frac{|\nabla %%@
u|^2}{|x|^{2m+2}}dx+\beta(W; R)\int_{B_{R}}W(x)\frac{|\nabla
u|^2}{|x|^{2m}}dx.
\end{equation}

 We also establish a more general version of equation (\ref{Rellich.1}). Assuming again that
$\frac{rW'(r)}{W(r)}$  decreases to $-\lambda$ on $(0, R)$, and provided  $m\leq \frac{n-4}{2}$ %%@
and $\frac{n}{2}+m\geq \lambda \geq n-2m-4$, we then have  for all
$u \in C^{\infty}_{0}(B_R)$,
\begin{eqnarray}\label{ex-gen-hr}
\int_{B_R}\frac{|\Delta u|^{2}}{|x|^{2m}}dx &\geq& \frac{(n+2m)^{2}(n-2m-4)^2}{16} \int_{B_R}\frac{u^2}{|x|^{2m+4}}dx %%@
\nonumber\\
&&\quad+\beta (W;  R)(\frac{(n+2m)^2}{4}+\frac{(n-2
m-\lambda-2)^2}{4}) \int_{B_R}\frac{W(x)}{|x|^{2m+2}}u^2 dx.
\end{eqnarray}

\section{General Hardy Inequalities}

Here is the main result of this section.
\begin{theorem} \label{main} Let $V$ and $W$ be positive radial $C^1$-functions on  %%@
$B_R\backslash \{0\}$, where $B_R$ is a ball centered at %%@
zero with radius $R$ ($0<R\leq +\infty$) in $\R^n$ ($n \geq 1$). Assume that $\int^{a}_{0}\frac{1}{r^{n-1}V(r)}dr=+\infty$ and %%@
$\int_{0}^{a}r^{n-1}V(r)dr<\infty$ for some $0<a<R$. Then the following two statements are equivalent:

\begin{enumerate}
\item The ordinary differential equation
\[
\hbox{ $({\rm B}_{V,W})$  \quad \quad \quad \quad \quad  \quad \quad \quad \quad \quad \quad \quad %%@
\quad \quad \quad $y''(r)+(\frac{n-1}{r}+\frac{V_r(r)}{V(r)})y'(r)+\frac{W(r)}{V(r)}y(r)=0$  \quad %%@
\quad \quad \quad \quad  \quad \quad \quad \quad \quad  \quad \quad \quad \quad \quad \quad \quad %%@
\quad \quad \quad}
\]
 has a positive solution on the interval $(0, R]$ with $\theta:=V(R)\frac{\phi'(R)}{\phi(R)}$.

\item For all $u \in H^{1}(B_R)$
\begin{equation*}\label{2dim-in}
\hbox{ $({\rm H}_{V,W})$ \quad \quad \quad \quad \quad \quad \quad \quad \quad \quad \quad \quad %%@
\quad   $\int_{B_R}V(x)|\nabla u(x) |^{2}dx \geq \int_{B_R} W(x)u^2dx$+$\theta \int_{\partial B} u^2$ ds.\quad \quad \quad \quad \quad %%@
\quad \quad  \quad \quad \quad \quad \quad \quad \quad \quad \quad
\quad}
\end{equation*}
\end{enumerate}

\end{theorem}

%\begin{theorem} \label{main-Rn} Let $V$ and $W$ be positive radial $C^1$-functions on a ball %%@
%$\R^n\backslash \{0\}$, with $n \geq 1$. Assume that $\int^{a}_{0}\frac{1}{r^{n-1}V(r)}dr=+\infty$ %%@
%and $\int_{0}^{a}r^{n-1}V(r)dr<\infty$ for some $a>0$, then the following two statements are %%@
%equivalent:
%\begin{enumerate}
%\item The ordinary differential equation
%\[
%\hbox{ $({\rm B'}_{V,W})$  \quad \quad \quad \quad \quad  \quad \quad \quad \quad \quad \quad %%@
%\quad \quad \quad \quad $y''(r)+(\frac{n-1}{r}+\frac{V_r(r)}{V(r)})y'(r)+\frac{W(r)}{V(r)}y(r)=0$  %%@
%\quad \quad \quad \quad \quad  \quad \quad \quad \quad \quad  \quad \quad \quad \quad \quad \quad %%@
%\quad \quad \quad \quad}
%\]
% has a positive solution on the interval $(0, \infty)$.
 %\item For all $u \in C_{0}^{\infty}(\R^n)$
%\begin{equation*}\label{2dim-in}
%\hbox{ $({\rm H'}_{V,W})$ \quad \quad \quad \quad \quad \quad \quad \quad \quad \quad \quad \quad %%@
%\quad   $\int_{\R^n}V(x)|\nabla u(x) |^{2}dx \geq \int_{\R^n} W(x)u^2dx$.\quad \quad \quad \quad %%@
%\quad \quad \quad  \quad \quad \quad \quad \quad \quad \quad \quad \quad \quad}
%\end{equation*}
%\end{enumerate}
%\end{theorem}

The above theorem allows to generalize all Hardy type inequalities
on $H^{1}_{0}(\Omega)$ to a corresponding inequality on
$H^{1}(\Omega)$. For instance we can get the following general form
of the Caffarelli-Kohn-Nirenberg inequalities.
\begin{corollary}
Assume $B$ is the ball of radius $R$ and and centered at zero in
$\R^n$. If $a\leq \frac{n-2}{2}$, then
\begin{equation}\label{no-improve}
\int_{B}|x|^{-2a}|\nabla u(x) |^{2}dx \geq (\frac{n-2a-2}{2})^2\int_{B}|x|^{-2a-2}u^2 %%@
dx-\frac{(n-2a-2)R^{-2a-1}}{2}\int_{\partial B} u^2dx,
\end{equation}
\end{corollary}
for all $u \in H^{1}(B)$.\\

To prove Theorem \ref{main} we shall need the following lemma.
\begin{lemma} \label{super} Let $V$ and $W$ be positive radial $C^1$-functions  on a ball %%@
$B\backslash \{0\}$, where $B$ is a ball with radius $R$ in $\R^n$ ($n \geq 1$) and centered at %%@
zero. Assume
\begin{eqnarray*}
\hbox{$\int_{B}\left(V(x)|\nabla u|^{2}-W(x)|u|^{2}\right)dx-\theta \int_{\partial B}u^2 ds\geq 0$ for all $u \in %%@
H^{1}(B)$, }
\end{eqnarray*}
 for some $\theta<0$. Then there exists a $C^{2}$-supersolution to the following linear elliptic equation
\begin{eqnarray}\label{pde}
-{\rm div}(V(x)\nabla u)-W(x)u&=&0, \ \ \ \ {\rm in} \ \  B, \\
u&>&0 \ \ \quad {\rm in} \ \  B \setminus \{0\}, \\
 V\nabla u.\nu&=&\theta u \ \quad {\rm in} \ \  \partial B.
\end{eqnarray}
 \end{lemma}
{\bf Proof:} Define
 \begin{eqnarray*}
\lambda_{1}(V):=\inf \{\frac{\int_{B}V(x)|\nabla\psi|^{2}- %%@
W(x)|\psi|^{2}-\theta \int_{\partial B}u^2}{\int_{B}|\psi|^{2}}; \ \
\psi \in C^{\infty}_{0}(B \setminus \{0\}) \}.
\end{eqnarray*}
By our assumption $\lambda_1(V)\geq 0$. Let $(\phi_{n}, \lambda^{n}_{1})$ be the first eigenpair %%@
for the problem
\begin{eqnarray*}
(L-\lambda_{1}(V)-\lambda^{n}_{1})\phi_{n}&=&0 \ \ on \ \ B \setminus B_{\frac{R}{n}}\\
\phi_n&=&0 \ \ on\ \ \partial B_{\frac{R}{n}} \\
V\nabla \phi_{n}.\nu &=&\theta \phi_{n} \ \ on\ \ \partial B,
\end{eqnarray*}
where $Lu=-{\rm div}(V(x)\nabla u)- W(x) u$, and $B_{\frac{R}{n}}$ is a ball of radius %%@
$\frac{R}{n}$, $n\geq 2$ . The eigenfunctions can be chosen in such a way that $\phi_{n}>0$ on %%@
$B \setminus B_{\frac{R}{n}}$ and $\varphi_{n}(b)=1$, for some $b \in B$ with %%@
$\frac{R}{2}<|b|<R$.

Note that $\lambda^{n}_{1}\downarrow 0$ as $n \rightarrow \infty$. Harnak's inequality yields that %%@
for any compact subset $K$, $\frac{{\rm max}_{K}\phi_{n}}{{\rm min}_{K}\phi_{n}}\leq C(K)$ with %%@
the later constant being independent of $\phi_{n}$.  Also standard elliptic estimates also yields %%@
that the family $(\phi_{n})$ have also uniformly bounded derivatives on the compact sets %%@
$B-B_{\frac{R}{n}}$.  \\
Therefore, there exists a subsequence $(\varphi_{n_{l_{2}}})_{l_{2}}$ of ($\varphi_{n})_{n}$ such %%@
that $(\varphi_{n_{l_{2}}})_{l_{2}}$ converges to some $\varphi_{2} \in C^{2}(B \setminus %%@
B(\frac{R}{2}))$.  Now consider  $(\varphi_{n_{l_{2}}})_{l_{2}}$ on $B \setminus %%@
B(\frac{R}{3})$.  Again there exists a subsequence  $(\varphi_{n_{l_{3}}})_{l_{3}}$ of %%@
$(\varphi_{n_{l_{2}}})_{l_{2}}$ which converges to $\varphi_{3} \in C^{2}( B \setminus %%@
B(\frac{R}{3}))$,  and $\varphi_{3}(x)=\varphi_{2}(x)$ for all $x \in B \setminus %%@
B(\frac{R}{2})$. By repeating this argument we get a supersolution $\varphi \in C^{2}( B %%@
\setminus\{ 0\})$ i.e. $L\varphi \geq 0$, such that $\varphi>0$ on $B \setminus \{0\}$ and $V\nabla \phi.\nu =\theta \phi$  on $\partial B$. %%@
\hfill $\square$\\

\noindent {\bf Proof of Theorem \ref{main}:} First we prove that 1) implies 2). Let $\phi \in %%@
C^{1}(0,R]$ be a solution of $(B_{V,W})$   such that $\phi (x)>0$ for all $x \in (0,R)$. Define %%@
$\frac{u(x)}{\varphi(|x|)}= \psi(x)$. Then
\[ |\nabla u|^2=(\varphi'(|x|))^2 \psi^2(x)+2\varphi'(|x|)\varphi(|x|)\psi(x)\frac{x}{|x|}.\nabla %%@
\psi+\varphi^2(|x|)|\nabla \psi|^2.\]
Hence,
\[ V(|x|)|\nabla u|^2\geq V(|x|)(\varphi'(|x|))^2 %%@
\psi^2(x)+2V(|x|)\varphi'(|x|)\varphi(|x|)\psi(x)\frac{x}{|x|}.\nabla \psi(x).\]
Thus, we have
\[\int_{B} V(|x|)|\nabla u|^2 dx \geq \int_{B} V(|x|) (\varphi'(|x|))^2 \psi^2(x)dx+ \int_{B}2 %%@
V(|x|) \varphi'(|x|)\varphi(|x|)\psi(x) \frac{x}{|x|}.\nabla \psi dx.\]
Let $B_{\epsilon}$ be a ball of radius $\epsilon$ centered at the origin. Integrate by parts to %%@
get
\begin{eqnarray*}
\int_{B} V(|x|)|\nabla u|^2 dx &\geq& \int_{B} V(|x|) (\varphi'(|x|))^2 %%@
\psi^2(x)dx+\int_{B_{\epsilon}}2V(|x|) \varphi'(|x|)\varphi(|x|)\psi(x) \frac{x}{|x|}.\nabla \psi %%@
dx\\
&+&\int_{B\backslash B_{\epsilon}}2 V(|x|) \varphi'(|x|)\varphi(|x|)\psi(x) \frac{x}{|x|}.\nabla %%@
\psi dx\\
&=&\int_{B_{\epsilon}} V(|x|) (\varphi'(|x|))^2 \psi^2(x)dx+\int_{B_{\epsilon}}2V(|x|) %%@
\varphi'(|x|)\varphi(|x|)\psi(x) \frac{x}{|x|}.\nabla \psi dx\\
&-&\int_{B\backslash B_{\epsilon}}\left\{\big(V(|x|) %%@
\varphi''(|x|)\varphi(|x|)+(\frac{(n-1)V(|x|)}{r}+V_r(|x|))\varphi'(|x|)\varphi(|x|)\big)\psi^2(x)
\right\}dx\\
&+&\int_{\partial (B \backslash B_{\epsilon})}V(|x|)\varphi'(|x|)\varphi(|x|)\psi^2(x)ds
\end{eqnarray*}
Let $\epsilon \rightarrow 0$ and use Lemma 2.3 in \cite{GM2} and the fact that $\phi$ is a solution %%@
of   $(D_{v,w})$ to get
\begin{eqnarray*}
\int_{B} V(|x|)|\nabla u|^2 dx &\geq& %%@
-\int_{B}[V(|x|)\varphi''(|x|)+(\frac{(n-1)V(|x|)}{r}+V_r(|x|))\varphi'(|x|)]\frac{u^2(x)}{
\varphi(|x|)}dx\\
&=&\int_{B}W(|x|)u^2(x)dx-\theta\int_{\partial B}u^2 ds.
\end{eqnarray*}
To show that 2) implies 1), we assume that inequality (${\rm H}_{V,W}$) holds on a ball $B$ of %%@
radius $R$, and then apply Lemma \ref{super} to obtain a $C^{2}$-supersolution for the equation %%@
(\ref{pde}). Now take the surface average of $u$, that is
\begin{equation}\label{sup-eq}
y(r)=\frac{1}{n\omega_{w} r^{n-1}}\int_{\partial B_{r}} u(x)dS=\frac{1}{n\omega_{n}} %%@
\int_{|\omega|=1}u(r\omega)d\omega >0,
\end{equation}
where $\omega_{n}$ denotes the volume of the unit ball in $R^{n}$. We may assume that the unit %%@
ball is contained in $B$ (otherwise we just use a smaller ball). It
is easy to see that $V(R)\frac{y'(R)}{y(R)}=\theta$. We clearly have
\begin{equation}
y''(r)+\frac{n-1}{r}y'(r)= \frac{1}{n\omega_{n}r^{n-1}}\int_{\partial B_{r}}\Delta u(x)dS.
\end{equation}
Since $u(x)$  is a supersolution of (\ref{pde}), we have
\[\int_{\partial B_{r}}div(V(|x|)\nabla u)ds-\int_{\partial B}W(|x|)udx\geq 0,\]
and therefore,
\[V(r)\int_{\partial B_{r}}\Delta u dS -V_r(r)\int_{\partial B_{r}} \nabla u.x %%@
ds-W(r)\int_{\partial B_{r}}u(x)ds\geq 0.\]
It follows that
\begin{equation}
V(r)\int_{\partial B_{r}}\Delta u dS -V_r(r)y'(r)-W(r)y(r)\geq 0,
\end{equation}
and in view of (\ref{sup-eq}),  we see that $y$ satisfies the inequality
\begin{equation}\label{ode}
V(r)y''(r)+(\frac{(n-1)V(r)}{r}+V_r(r))y'(r)\leq -W(r)y(r), \ \ \ \ for\ \ \ 0<r<R,
\end{equation}
that is it is a positive supersolution $y$ for $(B_{V,W})$ with
$V(R)\frac{y'(R)}{y(R)}=\theta$. Standard results in ODE now allow us to conclude that $(B_{V,W})$ has actually a positive solution %%@
on $(0, R)$, and the proof of theorem \ref{main} is now complete.
\hfill $\Box$\\

An immediate application of Theorem 2.6 in \cite{GM2} and Theorem
\ref{main} is the following very general Hardy inequality.

\begin{theorem} \label{super.hardy} Let $V(x)=V(|x|)$ be a strictly positive radial function on a smooth domain $\Omega$ containing %%@
$0$ such that  $R=\sup_{x \in \Omega} |x|$. Assume
%Let $V$ be an strictly positive $C^1$-function on $(0,R)$ such %%@
that for some $\lambda \in \R$
\begin{equation}
\hbox{$\frac{rV_r(r)}{V(r)}+\lambda \geq 0$ on $(0, R)$ and $\lim\limits_{r\to %%@
0}\frac{rV_r(r)}{V(r)}+\lambda =0$.}
\end{equation}

If $\lambda \leq n-2$, then  the following inequality holds  for any
Bessel potential $W$ on $(0, R)$:
\begin{eqnarray*}\label{v-hardy}
\int_{\Omega}V(x)|\nabla u|^{2}dx &\geq&
(\frac{n-\lambda-2}{2})^2\int_{\Omega}\frac{V(x)}{|x|^2}u^{2}dx +\beta (W; R)\int_{\Omega} %%@
V(x)W(x)u^{2}dx
\\ &+& V(R)(\frac{\phi'(R)}{\phi(R)}-\frac{n-\lambda-2}{2R})
\int_{\partial B} u^2 \quad \ \ for \ \ u \in H^{1}(\Omega),
\end{eqnarray*}
where $\phi$ is the corresponding solution of $(B_{1,W})$.

\end{theorem}
{\bf Proof:} Under our assumptions, it is easy to see that
$y=r^{\frac{n-\lambda-2}{2}}\phi(r)$ is a positive super-solution of
$B_{(V,V(\frac{n-\lambda-2}{2})^2r^{-2}+W)}$. Now apply  Theorem 2.6
in \cite{GM2} and Theorem \ref{main} to complete the proof. \hfill
$\Box$

\section{General Hardy-Rellich inequalities}
Let $0 \in \Omega \subset R^n$ be a smooth domain, and denote
\[C^{k}_{r}(\bar{\Omega})=\{v \in C^{k}(\bar{\Omega}): \mbox{v is radial }\}.\]

We start by considering  a general inequality for radial functions.

\begin{theorem} \label{mainrad.hr} Let $V$ and $W$ be positive radial $C^1$-functions on a ball %%@
$B\backslash \{0\}$, where $B$ is a ball with radius $R$ in $\R^n$ ($n \geq 1$) and centered at %%@
zero. Assume $\int^{R}_{0}\frac{1}{r^{n-1}V(r)}dr=\infty$ and $\lim_{r \rightarrow %%@
0}r^{\alpha}V(r)=0$ for some $\alpha< n-2$. Then the following
statements are equivalent:
\begin{enumerate}

\item $(V, W)$ is a Bessel pair on $(0, R)$ with $\theta:= V(R)\frac{\phi'(R)}{\phi(R)}$, where $\phi$ is the corresponding solution of
$(B_{(V,W)})$.

\item If $\lim_{r \rightarrow 0}r^{\alpha}V(r)=0$ for some $\alpha< n-2$, then the above are %%@
equivalent to
\[
\hbox{$\int_{B}V(x)|\Delta u |^{2}dx \geq  \int_{B} W(x)|\nabla  %%@
u|^{2}dx+(n-1)\int_{B}(\frac{V(x)}{|x|^2}-\frac{V_r(|x|)}{|x|})|\nabla
u|^2dx+(\theta+(n-1)V(R))\int_{\partial B}|\nabla u|^2$,}
\]
for all radial $u \in C^{\infty}(\bar{B})$.

\end{enumerate}

\end{theorem}
{\bf Proof:}  Assume $u \in C^{\infty}_{r}(\bar{B})$ and observe
that
\[\int_{B}V(x)|\Delta u %%@
|^{2}dx=n\omega_{n}\{\int^{R}_{0}V(r)u_{rr}^{2}r^{n-1}dr+(n-1)^2\int^{R}_{0}V(r)\frac{u^{2}_{r}}{r
^{2 }}r^{n-1}dr +2(n-1)\int^{R}_{0}V(r)uu_rr^{n-2}dr\}.\] Setting
$\nu=u_{r}$, we then have
\[\int_{B}V(x)|\Delta u |^{2}dx=\int_{B}V(x)|\nabla \nu |^{2}dx+(n-1) %%@
\int_{B}(\frac{V(|x|)}{|x|^2}-\frac{V_r(|x|)}{|x|})|\nu|^{2}dx+(n-1)V(R)\int_{\partial
B}|\nu|^2ds.
\] Thus, $({\rm HR}_{V,W})$ for radial functions is equivalent to
\[\int_{B}V(x)|\nabla \nu |^{2}dx\geq \int_{B}W(x)\nu^2 dx.\]

It therefore follows from Theorem \ref{main} that 1) and 2) are  equivalent. \hfill $\Box$\\

\subsection{The non-radial case}
The decomposition of a function into its spherical harmonics will be
one of our tools to prove our results. This idea has also been used
in \cite{TZ} and \cite{GM2}. Let $u \in
C^{\infty}(\bar{B})$. By decomposing $u$ into spherical %%@
harmonics we get
\[
\hbox{$u=\Sigma^{\infty}_{k=0}u_{k}$ where
$u_{k}=f_{k}(|x|)\varphi_{k}(x)$}
\]
 and $(\varphi_k(x))_k$ are the orthonormal eigenfunctions of the Laplace-Beltrami operator  %%@
with corresponding eigenvalues $c_{k}=k(N+k-2)$, $k\geq 0$. The functions $f_{k}$ belong to %%@
$u \in C^{\infty}([0,R])$, $f_{k}(R)=0$, and satisfy $f_{k}(r)=O(r^k)$ and $f'(r)=O(r^{k-1})$ as $r \rightarrow %%@
0$. In particular,
\begin{equation}\label{zero}
\hbox{ $\varphi_{0}=1$ and  $f_{0}=\frac{1}{n \omega_{n}r^{n-1}}\int_{\partial B_{r}}u ds=
\frac{1}{n \omega_{n}}\int_{|x|=1}u(rx)ds.$}
\end{equation}
We also have  for any $k\geq 0$, and any continuous real valued functions $v$ and $w$ on %%@
$(0,\infty)$,
\begin{equation}
\int_{R^n}V(|x|)|\Delta u_{k}|^{2}dx=\int_{R^n}V(|x|)\big( \Delta %%@
f_{k}(|x|)-c_{k}\frac{f_{k}(|x|)}{|x|^2}\big)^{2}dx,
\end{equation}
and
\begin{equation}
\int_{R^n}W(|x|)|\nabla u_{k}|^{2}dx=\int_{R^n}W(|x|)|\nabla %%@
f_{k}|^{2}dx+c_{k}\int_{R^n}W(|x|)|x|^{-2}f^{2}_{k}dx.
\end{equation}

\begin{theorem} \label{main.hr} Let $V$ and $W$ be positive radial $C^1$-functions  on a ball %%@
$B\backslash \{0\}$, where $B$ is a ball with radius $R$ in $\R^n$ ($n \geq 1$) and centered at %%@
zero. Assume $\int^{R}_{0}\frac{1}{r^{n-1}V(r)}dr=\infty$ and $\lim_{r \rightarrow %%@
0}r^{\alpha}V(r)=0$ for some $\alpha<(n-2)$. If
\begin{equation}\label{main.con}
W(r)-\frac{2V(r)}{r^2}+\frac{2V_r(r)}{r}-V_{rr}(r)\geq 0 \ \ for \ \
0\leq r \leq R,
\end{equation}
and the ordinary differential equation $(B_{V,W})$ has a positive
solution $\phi$ on the interval $(0, R]$ such that
\begin{equation}\label{bd.con}
(n-1+R\frac{\phi'(R)}{\phi(R)})V(R)\geq 0,
\end{equation}
then the  following inequality holds
 for all $u \in H^{2}(B)$.
\begin{equation*}\label{gen-hardy}
\hbox{   $({\rm HR}_{V,W})$ \quad \quad  $\int_{B}V(x)|\Delta u |^{2}dx \geq  \int_{B} W(x)|\nabla %%@
u|^{2}dx+(n-1)\int_{B}(\frac{V(x)}{|x|^2}-\frac{V_r(|x|)}{|x|})|\nabla u|^2dx.$  \quad \quad \quad %%@
\quad \quad \quad \quad  }
\end{equation*}
Moreover, if  $\beta (V, W; R)\geq 1$, then the best constant is
given by
\begin{equation}
\hbox{$\beta (V, W; R)=\sup\big\{c; \, \,  ({\rm HR}_{V, cW})$
holds$\big\}$.}
\end{equation}
\end{theorem}

{\bf Proof:}  Assume that the equation  $(B_{V,W})$ has a positive solution on %%@
$(0,R]$. We prove that the inequality $(HR_{V,W})$ holds for all $u \in C^{\infty}_{0}(B)$ by %%@
frequently using that

\begin{equation}\label{1-dim}
\hbox{$\int^{R}_{0}V(r)|x'(r)|^2r^{n-1}dr \geq \int^{R}_{0}W(r)x^{2}(r)r^{n-1}dr+V(R)\frac{\phi'(R)}{\phi(R)}R^{n-1} (x(R))^2$ for all $x\in %%@
C^1(0, R]$.}
\end{equation}

Indeed,  for all $n\geq 1$ and $k\geq 0$ we have
\begin{eqnarray*}
\frac{1}{nw_n}\int_{R^n}V(x)|\Delta u_{k}|^{2}dx&=&\frac{1}{nw_n}\int_{R^n}V(x)\big( \Delta %%@
f_{k}(|x|)-c_{k}\frac{f_{k}(|x|)}{|x|^2}\big)^{2}dx\\
&=&
\int^{R}_{0}V(r)\big(f_{k}''(r)+\frac{n-1}{r}f_{k}'(r)-c_{k}\frac{f_{k}(r)}{r^2}\big)^{2}r^{n-1}dr
\\
&=&\int^{R}_{0}V(r)(f_{k}''(r))^{2}r^{n-1}dr+(n-1)^{2} \int^{R}_{0}V(r)(f_{k}'(r))^{2}r^{n-3}dr\\
&&+c^{2}_{k} \int^{R}_{0}V(r)f_{k}^{2}(r)r^{n-5}
+ 2(n-1) \int^{R}_{0}V(r)f_{k}''(r)f_{k}'(r)r^{n-2}\\
&&-2c_{k}
\int^{R}_{0}V(r)f_{k}''(r)f_{k}(r)r^{n-3}dr
- 2c_{k}(n-1) \int^{R}_{0}V(r)f_{k}'(r)f_{k}(r)r^{n-4}dr.
\end{eqnarray*}
Integrate by parts and use (\ref{zero}) for $k=0$ to get
\begin{eqnarray}
\frac{1}{n\omega_{n}}\int_{R^n}V(x)|\Delta u_{k}|^{2}dx&=&  %%@
\int^{R}_{0}V(r)(f_{k}''(r))^{2}r^{n-1}dr+(n-1+2c_{k}) \int^{R}_{0}V(r)(f_{k}'(r))^{2}r^{n-3}dr %%@
\label{piece.0}\\
&+&
(2c_{k}(n-4)+c^{2}_{k})\int^{R}_{0}V(r)r^{n-5}f_{k}^{2}(r)dr-(n-1)\int^{R}_{0}V_r(r)r^{n
-2}(f_{k}')^{2}(r)dr\nonumber\\
&-&c_{k}(n-5)\int^{R}_{0}V_r(r)f_{k}^2(r)r^{n-4}dr-c_{k}\int^{R}_{0}V_{rr}(r)f_{k}^2(r)r^{n-3}dr.\\
&+&(n-1)V(R)(f_{k}'(R))^2R^{n-2} \nonumber
\end{eqnarray}
Now define $g_{k}(r)=\frac{f_{k}(r)}{r}$ and note that $g_{k}(r)=O(r^{k-1})$ for all $k\geq 1$. We %%@
have
\begin{eqnarray*}
\int^{R}_{0}V(r)(f_{k}'(r))^{2}r^{n-3}&=&\int^{R}_{0}V(r)(g_{k}'(r))^{2}r^{n-1}dr+\int^{R}_{0}2V(r
)g
_{k}(r)g_{k}'(r)r^{n-2}dr+\int^{R}_{0}V(r)g_{k}^{2}(r)r^{n-3}dr\\
&=&\int^{R}_{0}V(r)(g_{k}'(r))^{2}r^{n-1}dr-(n-3)\int^{R}_{0}V(r)g_{k}^{2}(r)r^{n-3}dr
-\int_{0}^{R}V_r(r)g^2_{k}(r)r^{n-2}dr\\
\end{eqnarray*}
Thus,
\begin{equation}\label{g1}
\int^{R}_{0}V(r)(f_{k}'(r))^{2}r^{n-3}\geq %%@
\int^{R}_{0}W(r)f_{k}^{2}(r)r^{n-3}dr-(n-3)\int^{R}_{0}V(r)f_{k}^{2}(r)r^{n-5}dr-\int_{0}^{R}V_r(r%%@
%%@
)f^2_{k}(r)r^{n-4}dr.
\end{equation}
Substituting $2c_k\int^{R}_{0}V(r)(f_{k}'(r))^{2}r^{n-3}$ in (\ref{piece.0}) by its lower estimate %%@
in the last inequality (\ref{g1}), we get
\begin{eqnarray*}
\frac{1}{n\omega_{n}}\int_{R^n}V(x)|\Delta u_{k}|^{2}dx&\geq&
\int^{R}_{0}W(r)(f_{k}'(r))^{2}r^{n-1}dr+\int^{R}_{0}W(r)(f_{k}(r))^{2}r^{n-3}dr \\  &+&(n-1) %%@
\int^{R}_{0}V(r)(f_{k}'(r))^{2}r^{n-3}dr+c_{k}(n-1) \int^{R}_{0}V(r)(f_{k}(r))^{2}r^{n-5}dr\\
&-&(n-1)\int^{R}_{0}V_r(r)r^{n-2}(f_{k}')^{2}(r)dr-c_{k}(n-1)\int^{R}_{0}V_r(r
)r^{n
-4}(f_{k})^{2}(r)dr\\
&+&c_{k}(c_{k}-(n-1))\int^{R}_{0}V(r)r^{n-5}f_{k}^{2}(r)dr\\
&+&c_{k}\int^{R}_{0}(W(r)-\frac{2V(r)}{r^2}+\frac{2V_r(r)}{r}-V_{rr}(r))f^2_{k}(r)r^{n-3}dr\\
&+&(n-1)V(R)(f_{k}'(R))^2R^{n-2}+V(R)\frac{\phi'(R)}{\phi(R)}R^{n-1}(f_{k}'(R))^2\\
\end{eqnarray*}
The proof is now complete since the last two terms are non-negative
by our assumptions. \hfill $\square$
\begin{remark}\rm
In order to apply the above theorem to
\[
\hbox{$V(x)=|x|^{-2m}$}
\]
we see that even in the simplest case $V\equiv 1$ condition (\ref{main.con}) reduces to %%@
$(\frac{n-2}{2})^2|x|^{-2}\geq 2|x|^{-2}$, which is then guaranteed
only if $n\geq 5$. More generally,  if $V(x)=|x|^{-2m}$, then in
order to satisfy (\ref{main.con}) we need to have
\begin{equation}\label{restrict}
\frac{-(n+4)-2\sqrt{n^2-n+1}}{6} \leq m\leq
\frac{-(n+4)+2\sqrt{n^2-n+1}}{6}.
\end{equation}
Also to satisfy the condition $(\ref{bd.con})$ we need to have
$m>-\frac{n}{2}$. Thus for $m$ satisfying $(\ref{restrict})$ the
inequality
\begin{equation}\label{gm-hr.0}
\int_{B_{R}}\frac{|\Delta u|^2}{|x|^{2m}}\geq  (\frac{n+2m}{2})^2\int_{B_{R}}\frac{|\nabla %%@
u|^2}{|x|^{2m+2}}dx.
\end{equation}
for all $u \in H^{2}(B_{R})$. Moreover, $(\frac{n+2m}{2})^2$ is the
best constant. We shall see however that %%@
this inequality remains true without condition (\ref{restrict}), but with a  constant
that is sometimes different from $(\frac{n+2m}{2})^2$ in the cases where (\ref{restrict}) is not %%@
valid. For example, if $m=0$, then the best constant is $3$ in dimension $4$ and $\frac{25}{36}$ %%@
in dimension $3$.
\end{remark}
%\begin{remark} \rm In view of the above remarks, Theorem \ref{main.hr} in the case where $m=0$ and %%@
%$n\geq 5$ already includes Theorem 1.5 in \cite{TZ} as a special case. Moreover,  when %%@
%$V(x)=|x|^{-2m}$,  the above shows that Theorem 1.8 in \cite{TZ} is still valid if one replaces %%@
%the condition
%\[0\leq m\leq \frac{-(n+4)+2\sqrt{n^2-n+1}}{6},\]
%by
%\[\frac{-(n+4)-2\sqrt{n^2-n+1}}{6} \leq m\leq \frac{-(n+4)+2\sqrt{n^2-n+1}}{6}.
%\]
%\end{remark}

\subsection{The case of  power potentials $|x|^m$}

The general Theorem \ref{main.hr} allowed us to deduce inequality (\ref{gm-hr}) below for a %%@
restricted interval of powers $m$. We shall now prove that the same holds for all  %%@
$-\frac{n}{2}\leq m<\frac{n-2}{2}$. We start with the following
result.

\begin{theorem}
Assume $-\frac{n}{2}\leq m<\frac{n-2}{2}$ and $\Omega$ be a smooth
domain in $\R^n$, $n\geq 1$. Then
\begin{equation}
a_{n,m}=\inf \left\{\frac{\int_{B_{R}}\frac{|\Delta u|^2}{|x|^{2m}}dx}{\int_{B_{R}}\frac{|\nabla %%@
u|^2}{|x|^{2m+2}}dx};\, H^{2}(\Omega)\setminus
\{0\}\right\}=\inf \left\{\frac{\int_{B_{R}}\frac{|\Delta u|^2}{|x|^{2m}}dx}{\int_{B_{R}}\frac{|\nabla %%@
u|^2}{|x|^{2m+2}}dx};\, u\in H^{2}_{0}(\Omega)\setminus
\{0\}\right\}
\end{equation}

\end{theorem}
{\bf Proof.} Decomposing again $u \in C^{\infty}(\bar{B}_{R})$ into spherical harmonics; $u=\Sigma^{\infty}_{k=0}u_{k}$, %%@
where $u_{k}=f_{k}(|x|)\varphi_{k}(x)$, one has
\begin{eqnarray}\label{N1}
\int_{\R^n}\frac{|\Delta %%@
u_{k}|^2}{|x|^{2m}}dx&=&\int_{\R^n}|x|^{-2m}(f''_{k}(|x|))^2dx+\left((n-1)(2m+1)+2c_{k}\right) %%@
\int_{\R^n}|x|^{-2m-2}(f_{k}')^2dx\\
&+&c_{k}(c_{k}+(n-4-2m)(2m+2))\int_{\R^n}|x|^{-2m-4}(f_{k})^2dx+(n-1)R^{n-2m-2}(f_{k}'(R))^2,
\nonumber
\end{eqnarray}
and
\begin{equation}\label{N2}
\int_{\R^n}\frac{|\nabla u_{k}|^2}{|x|^{2m+2}}dx=\int_{\R^n}|x|^{-2m-2}(f_{k}')^2 %%@
dx+c_{k}\int_{\R^n}|x|^{-2m-4}(f_{k})^2dx.
\end{equation}
The rest of the proof follows from the inequality
$(\ref{no-improve})$ and an argument similar to that of Theorem 6.1
in \cite{GM2}. \hfill $\Box$

\begin{remark}
The constant $a_{n,m}$ has been computed explicitly in \cite{GM2}
(Theorem 6.1).
\end{remark}

 \begin{theorem}\label{gm.hr}
Suppose $n\geq 1$ and $-\frac{n}{2}\leq m<\frac{n-2}{2}$, and  $W$ is a Bessel potential on $B_{R} \subset R^n$ with $n\geq 3$ and $\phi$ is the corresponding solution for the $(B_{1,W})$.  %%@
If
\[R\frac{\phi'(R)}{\phi(R)}\geq -\frac{n}{2}-m,\]
then for all $u \in H^{2}(B_{R})$ we
have
\begin{equation}\label{gm-hr}
\int_{B_{R}}\frac{|\Delta u|^2}{|x|^{2m}}\geq a_{n,m}\int_{B_{R}}\frac{|\nabla %%@
u|^2}{|x|^{2m+2}}dx+\beta(W;  R)\int_{B_{R}}W(x)\frac{|\nabla u|^2}{|x|^{2m}}dx,
\end{equation}
where
\[
a_{n,m}=\inf \left\{\frac{\int_{B_{R}}\frac{|\Delta u|^2}{|x|^{2m}}dx}{\int_{B_{R}}\frac{|\nabla %%@
u|^2}{|x|^{2m+2}}dx};\, u\in H^{2}(B_{R}) \setminus \{0\}\right\}.
\]
Moreover $\beta(W;  R)$ and $a_{m,n}$ are the best constants to be
computed in the appendix.
\end{theorem}
{\bf Proof:}  Assuming  the inequality
\[\int_{B_{R}}\frac{|\Delta u|^2}{|x|^{2m}}\geq a_{n,m}\int_{B_{R}}\frac{|\nabla %%@
u|^2}{|x|^{2m+2}}dx,\]
holds for all $u \in C^{\infty}(\bar{B}_{R})$, we shall prove that it can be improved by any Bessel %%@
potential $W$. We will use the following inequality in the proof which follows directly %%@
from the inequality (\ref{no-improve}) with n=1.
\begin{equation}\label{freq-in}
\int_{0}^{R}r^{\alpha}(f'(r))^{2}dr\geq %%@
(\frac{\alpha-1}{2})^2\int_{0}^{R}r^{\alpha-2}f^2(r)dr+\beta(W;R)\int_{0}^{R}r^{\alpha}W(r)
f^2(r)dr+(\frac{\phi'(R)}{\phi(R)}-\frac{\alpha-1}{2R})R^{\alpha},
\end{equation}
for $\alpha \geq 1$ and for all $f \in C^{\infty}(0,R]$, where both $(\frac{\alpha-1}{2})^2$ and $\beta(W;R)$ are best %%@
constants. Decompose $u \in C^{\infty}(\bar{B}_{R})$ into its spherical harmonics $ %%@
\Sigma^{\infty}_{k=0}u_{k}$, where $u_{k}=f_{k}(|x|)\varphi_{k}(x)$. We evaluate %%@
$I_k=\frac{1}{nw_n}\int_{R^n} \frac{|\Delta u_{k}|^2}{|x|^{2m}}dx$ in the following way
\begin{eqnarray*}
 I_k&=&\int_{0}^{R}r^{n-2m-1}(f''_{k}(r))^2dr+[(n-1)(2m+1)+2c_{k}]\int^{R}_{0}r^{n-2m-3}(f_{k}')^2
dr\\
&&+c_{k}[c_{k}+(n-2m-4)(2m+2)]\int_{0}^{R}r^{n-2m-5}(f_{k}(r))^2dr\\
&+& (n-1)R^{n-2m-2}(f_{k}'(R))^2\\
&\geq& \beta %%@
(W)\int^{R}_{0}r^{n-2m-1}W(x)(f_{k}')^2dr+[(\frac{n+2m}{2})^2+2c_{k}]\int^{R}_{0}r^{n-2m-3}(f_{k}'%%@
)^2dr\\
&&+c_{k}[c_{k}+(n-2m-4)(2m+2)]\int_{0}^{R}r^{n-2m-5}(f_{k}(r))^2dr\\
&\geq& \beta (W)\int^{R}_{0}r^{n-2m-1}W(x)(f_{k}')^2dr+a_{n,m}\int^{R}_{0}r^{n-2m-3}(f_{k}')^2dr\\
&&+\beta (W)[(\frac{n+2m}{2})^2+2c_{k}-a_{n,m}]\int^{R}_{0}r^{n-2m-3}W(x)(f_{k})^2dr\\
&&+\big((\frac{n-2m-4}{2})^2[(\frac{n+2m}{2})^2+2c_{k}-a_{n,m}]+c_{k}[c_{k}+(n-2m-4)(2m+2)]\big)
\int_{0}^{R}r^{n-2m-5}(f_{k}(r))^2dr.\\
\end{eqnarray*}
Now by (115) in \cite{GM2} we have
\[\big((\frac{n-2m-4}{2})^2[(\frac{n+2m}{2})^2+2c_{k}-a_{n,m}]+c_{k}[c_{k}+(n-2m-4)(2m+2)]\geq %%@
c_{k}a_{n,m},\]
for all $k\geq 0$. Hence, we have
\begin{eqnarray*}
I_{k}&\geq&a_{n,m}\int^{R}_{0}r^{n-2m-3}(f_{k}')^2dr+a_{n,m}c_{k}\int_{0}^{R}r^{n-2m-5}(f_{k}(r))^%%@
2dr\\
&&+\beta (W)\int^{R}_{0}r^{n-2m-1}W(x)(f_{k}')^2dr+\beta %%@
(W)[(\frac{n+2m}{2})^2+2c_{k}-a_{n,m}]\int^{R}_{0}r
^{n-2m-3}W(x)(f_{k})^2dr\\
&\geq& a_{n,m}\int^{R}_{0}r^{n-2m-3}(f_{k}')^2dr+a_{n,m}c_{k} %%@
\int_{0}^{R}r^{n-2m-5}(f_{k}(r))^2dr\\
&&+ \beta (W)\int^{R}_{0}r^{n-2m-1}W(x)(f_{k}')^2dr+
\beta (W)c_{k}\int^{R}_{0}r^{n-2m-3}W(x)(f_{k})^2dr\\
&=&a_{n,m}\int_{B_{R}}\frac{|\nabla u|^2}{|x|^{2m+2}}dx+\beta (W)\int_{B_{R}}W(x)\frac{|\nabla %%@
u|^2}{|x|^{2m}}dx.
\end{eqnarray*}
\hfill $\Box$
\\
In the following theorem we prove a very general class of weighted
Hardy-Rellich inequalities on $H^{2}(\Omega)\cap H^{1}_{0}$.

\begin{theorem}\label{super.hardy-rellich}
Let $\Omega$ be a smooth domain in $R^{n}$ with $n\geq 1$ and let $V \in C^{2}(0,R=:\sup_{x \in %%@
\Omega}|x|)$ be a non-negative function that satisfies the following conditions:
\begin{equation}
\hbox{$V_r(r)\leq 0$\quad and \quad
$\int^{R}_{0}\frac{1}{r^{n-3}V(r)}dr=-\int^{R}_{0}\frac{1}{r^{n-4}V_r(r)}dr=+\infty$.}
\end{equation}
 There exists $\lambda_{1}, \lambda_{2} \in R$ such that
\begin{equation}
\hbox{$\frac{rV_r(r)}{V(r)}+\lambda_{1} \geq 0$ on $(0, R)$ and $\lim\limits_{r\to %%@
0}\frac{rV_r(r)}{V(r)}+\lambda_{1} =0$,}
\end{equation}
\begin{equation}
\hbox{$\frac{rV_{rr}(r)}{V_r(r)}+\lambda_{2} \geq 0$ on $(0, R)$ and $\lim\limits_{r\to %%@
0}\frac{rV_{rr}(r)}{V_r(r)}+\lambda_{2} =0$,}
\end{equation}
and
\begin{equation}\label{super.hr.con}
\hbox{$\left(\frac{1}{2}(n-\lambda_{1}-2)^2+3(n-3)\right)V(r)-(n-5)rV_r(r)-r^2V_{rr}(r)\geq 0$ for %%@
all $r \in (0,R)$. }
\end{equation}
If $\lambda_{1}\leq n$, then the following inequality holds:
\begin{eqnarray}\label{super.hr}
\int_{\Omega}V(|x|)|\Delta u|^2 dx&\geq& %%@
(\frac{(n-\lambda_{1}-2)^{2}}{4}+(n-1))\frac{(n-\lambda_{1}-4)^{2}}{4}\int_{\Omega}\frac{V(|x|)}{|%%@
x|^4}u^2 dx \nonumber\\
&&-\frac{(n-1)(n-\lambda_{2}-2)^{2}}{4}\int_{\Omega}\frac{V_r(|x|)}{|x|^3}u^2 dx.
\end{eqnarray}
\end{theorem}
{\bf Proof:} We have by Theorem \ref{super.hardy} and condition (\ref{super.hr.con}),
{\small
\begin{eqnarray*}
\frac{1}{n\omega_{n}}\int_{R^n}V(x)|\Delta u_{k}|^{2}dx&=&  %%@
\int^{R}_{0}V(r)(f_{k}''(r))^{2}r^{n-1}dr+(n-1+2c_{k}) \int^{R}_{0}V(r)(f_{k}'(r))^{2}r^{n-3}dr %%@
\label{piece}\\
&+&
(2c_{k}(n-4)+c^{2}_{k})\int^{R}_{0}V(r)r^{n-5}f_{k}^{2}(r)dr-(n-1)\int^{R}_{0}V_r(r)r^{n
-2}(f_{k}')^{2}(r)dr\\
&-&c_{k}(n-5)\int^{R}_{0}V_r(r)f_{k}^2(r)r^{n-4}dr-c_{k}\int^{R}_{0}V_{rr}(r)f_{k}^2(r)r^{n-3}dr\\
&+&(n-1)V(R)(f_{k}'(R))^2R^{n-2}\\
&\geq&\int^{R}_{0}V(r)(f_{k}''(r))^{2}r^{n-1}dr+(n-1) \int^{R}_{0}V(r)(f_{k}'(r))^{2}r^{n-3}dr %%@
\label{piece}\\
&-&
(n-1)\int^{R}_{0}V_r(r)r^{n-2}(f_{k}')^{2}(r)dr\\
&+&c_{k}\int_{0}^{R}\left(\left(\frac{1}{2}(n-\lambda_{1}-2)^2+3(n-3)\right)V(r)-(n-5)rV_r(r)-r^2V%%@
%%@
_{rr}(r) \right)f_{k}^2(r)r^{n-5}dr
%&\times& f_{k}^2(r)r^{n-5}dr,
\\
&+&(n-1)V(R)(f_{k}'(R))^2R^{n-2}
\end{eqnarray*}
}
The rest of the proof follows %%@
from the above inequality combined with Theorem \ref{super.hardy}. \hfill $\Box$
\begin{remark}\rm
Let $V(r)=r^{-2m}$ with $-\frac{n}{2}\leq m\leq \frac{n-4}{2}$. Then in order to satisfy condition %%@
(\ref{super.hr.con}) we must have $-1-\frac{\sqrt{1+(n-1)^2}}{2}\leq
m\leq \frac{n-4}{2}$. Since
$-1-\frac{\sqrt{1+(n-1)^2}}{2}\leq-\frac{n}{2}$,
 if $-\frac{n}{2}\leq m\leq \frac{n-4}{2}$ the inequality (\ref{super.hr}) gives the following weighted second order Rellich %%@
inequality:
\[\int_{B}\frac{|\Delta u|^2}{|x|^{2m}}dx\geq %%@
H_{n,m}\int_{B}\frac{u^2}{|x|^{2m+4}}dx \ \ u \in H^{2}(\Omega)\cap
H^{1}_{0}(\Omega),\] where
\begin{equation}
H_{n,m}:=(\frac{(n+2m)(n-4-2m)}{4})^2.
\end{equation}

\end{remark}

The following theorem includes a large class of improved
Hardy-Rellich inequalities as special cases.

\begin{theorem}\label{gen.hr} Let $-\frac{n}{2}\leq m\leq \frac{n-4}{2}$ and let $W(x)$ be a Bessel potential on a %%@
ball $B$ of radius $R$  in $R^n$ with radius $R$. Assume %%@
$\frac{W(r)}{W_r(r)}=-\frac{\lambda}{r}+f(r)$, where $f(r)\geq 0$ and $\lim_{r \rightarrow %%@
0}rf(r)=0$. If $\lambda \leq \frac{n}{2}+m$, then the following
inequality holds for all $u \in H^{2}\cap H^{1}_{0}(B)$
\begin{eqnarray}\label{ex-gen-hr}
\int_{B}\frac{|\Delta u|^{2}}{|x|^{2m}}dx &\geq& H_{n,m}\int_{B}\frac{u^2}{|x|^{2m+4}}dx %%@
\nonumber\\
&&\quad+\beta (W;  R)(\frac{(n+2m)^2}{4}+\frac{(n-2
m-\lambda-2)^2}{4})
\int_{B}\frac{W(x)}{|x|^{2m+2}}u^2 dx.
\end{eqnarray}
Moreover, both constants are the best constants.
\end{theorem}

{\bf Proof:} Again we will frequently use inequality (\ref{freq-in})  in the proof. Decomposing $u %%@
\in C^{\infty}(\bar{B}_{R})$ into spherical harmonics
$\Sigma^{\infty}_{k=0}u_{k}$, where
$u_{k}=f_{k}(|x|)\varphi_{k}(x)$, we can write
\begin{eqnarray*}
\frac{1}{n\omega_{n}}\int_{R^n} \frac{|\Delta %%@
u_{k}|^2}{|x|^{2m}}dx&=&\int_{0}^{R}r^{n-2m-1}(f''_{k}(r))^2dr+[(n-1)(2m+1)+2c_{k}]\int^{R}_{0}r^{%%@
n-2m-3}(f_{k}')^2dr\\
&&+c_{k}[c_{k}+(n-2m-4)(2m+2)]\int_{0}^{R}r^{n-2m-5}(f_{k}(r))^2dr\\
&+&(n-1)(f_{k}'(R))^2R^{n-2m-2}\\
&\geq&(\frac{n+2m}{2})^2\int^{R}_{0}r^{n-2m-3}(f_{k}')^2dr+ \beta (W;  %%@
R)\int^{R}_{0}r^{n-2m-1}W(x)(f_{k}')^2dr\\
&+&c_{k}[c_{k}+2(\frac{n-\lambda-4}{2})^2+(n-2m-4)(2m+2)]\int_{0}^{R}r^{n-2m-5}(f_{k}(r))^2dr\\
&+&(n-1)(f_{k}'(R))^2R^{n-2m-2},
\end{eqnarray*}
where we have used the fact that $c_k\geq 0$ to get the above inequality. We have
\begin{eqnarray*}
\frac{1}{n\omega_{n}}\int_{R^n} \frac{|\Delta u_{k}|^2}{|x|^{2m}}dx
&\geq& \beta_{n,m} \int^{R}_{0}r^{n-2m-5}(f_{k})^2dr\\
&&+\beta (W;  R)\frac{(n+2m)^{2}}{4}
\int^{R}_{0}r^{n-2m-3}W(x)(f_{k})^2dr\\
&& +\beta (W;  R)\int^{R}_{0}r^{n-2m-1}W(x)(f_{k}')^2dr\\
&\geq&\beta_{n,m}\int^{R}_{0}r^{n-2m-5}(f_{k})^2dr\\
&&+\beta (W;  R)(\frac{(n+2m)^2}{4}+\frac{(n-2
m-\lambda-2)^2}{4})
\int^{R}_{0}r^{n-2m-3}W(x)(f_{k})^2dr \\
&\geq& \frac{\beta_{n,m}}{n\omega_{n}}\int_{B}\frac{u_{k}^2}{|x|^{2m+4}}dx\\
&&+\frac{\beta (W;  R)}{n\omega_{n}}(\frac{(n+2m)^2}{4}+\frac{(n-2
m-\lambda-2)^2}{4})
\int_{B}\frac{W(x)}{|x|^{2m+2}}u_{k}^2 dx,
\end{eqnarray*}
by Theorem \ref{super.hardy}. Hence, (\ref{ex-gen-hr}) holds and the proof is complete. \hfill %%@
$\Box$\\

We shall now give a few immediate applications of the above in the
case where $m=0$ and $n\geq 3$.

\begin{theorem} \label{m=0}Assume $W$ is a Bessel potential on $B_{R} \subset R^n$ with $n\geq 3$ and $\phi$ is the corresponding solution for the $(B_{1,W})$.  %%@
If
\[R\frac{\phi'(R)}{\phi(R)}\geq -\frac{n}{2},\]
then for all $u \in H^{2}(B_{R})$ we have
\begin{equation}
\int_{B_{R}}|\Delta u|^2 dx\geq C(n)\int_{B_{R}}\frac{|\nabla
u|^2}{|x|^2}dx+\beta(W;R)\int_{B_{R}}W(x)|\nabla u|^2dx,
\end{equation}
where $C(3)=\frac{25}{36}$, $C(4)=3$ and $C(n)=\frac{n^2}{4}$ for
all $n\geq 5$. Moreover, $C(n)$ and $\beta(W;R)$ are best
constants.\\
\end{theorem}

\begin{corollary}
The following holds for any smooth bounded domain $\Omega$ in $R^n$ with %%@
$R=\sup_{x \in \Omega}|x|$, and any $u \in H^{2}(\Omega)$.
\begin{enumerate}
\item  Let $z_{0}$ be the first zero of the Bessel function $J_{0}(z)$ and  choose
$0<\mu<z_0$ so that

\begin{equation}\label{mu}
\mu\frac{J_{0}'(\mu)}{J_{0}(\mu)}=-\frac{n}{2}.
\end{equation}
Then
 \begin{equation}\label{gm-hardy-Rellich}
\hbox{$\int_{\Omega}|\Delta u |^{2}dx \geq C(n) \int_{\Omega}\frac{|\nabla %%@
u|^2}{|x|^{2}}dx+\frac{\mu^{2}}{R^2}\int_{\Omega}|\nabla u|^{2}dx$}
\end{equation}

\item For any $k\geq 1$, choose $\rho\geq R( e^{e^{e^{.^{.^{e(k-times)}}}}} )$ large enough so
that $R\frac{\phi'(R)}{\phi(R)}\geq -\frac{n}{2}$, where
\begin{equation}\label{log.def}
\phi=\big(
\prod^{j}_{i=1}log^{(i)}\frac{\rho}{|x|}\big)^{\frac{1}{2}}.
\end{equation}
Then we have
 \begin{equation}
 \int_{\Omega}|\Delta u(x) |^{2}dx \geq C(n)\int_{\Omega}\frac{|\nabla u|^2}{|x|^2} %%@
dx+\frac{1}{4}\sum^{k}_{j=1}\int_{\Omega}\frac{|\nabla u |^2}{|x|^2}\big( %%@
\prod^{j}_{i=1}log^{(i)}\frac{\rho}{|x|}\big)^{-2}dx,
 \end{equation}
 \item We have
 \begin{equation}
  \int_{\Omega}|\Delta u(x) |^{2}dx \geq C(n)\int_{\Omega}\frac{|\nabla u|^2}{|x|^2} %%@
dx+\frac{1}{4}\sum^{n}_{i=1}\int_{\Omega}\frac{|\nabla %%@
u|}{|x|^{2}}X^{2}_{1}(\frac{|x|}{R})X^{2}_{2
}(\frac{|x|}{R})...X^{2}_{i}(\frac{|x|}{R})dx.
\end{equation}

\end{enumerate}
\end{corollary}

The following is immediate from Theorem \ref{gen.hr} and from the fact that $\lambda=2$ for the %%@
Bessel potential under consideration.

\begin{corollary} Let $\Omega$ be a smooth bounded domain in $\R^n$, $n \geq 4$ and $R=\sup_{x \in %%@
\Omega}|x|$. Then the following holds for all $u \in H^{2}(\Omega)
\cap H^{1}_{0}(\Omega)$
\begin{enumerate}
\item  Choose $\rho\geq R( e^{e^{e^{.^{.^{e(k-times)}}}}} )$ so that $R\frac{\phi'(R)}{\phi(R)} \geq
-\frac{n}{2}$. Then
\begin{equation}
 \int_{\Omega}|\Delta u(x) |^{2}dx \geq \frac{n^2(n-4)^2}{16}\int_{\Omega}\frac{u^2}{|x|^4} %%@
dx+(1+\frac{n(n-4)}{8})\sum^{k}_{j=1}\int_{\Omega}\frac{u^2}{|x|^4}\big( %%@
\prod^{j}_{i=1}log^{(i)}\frac{\rho}{|x|}\big)^{-2}dx.
 \end{equation}
\item Let $X_{i}$ is defined as in the introduction, then
\begin{equation}
\int_{\Omega}|\Delta u(x) |^{2}dx \geq \frac{n^2(n-4)^2}{16}\int_{\Omega}\frac{u^2}{|x|^4} %%@
dx+(1+\frac{n(n-4)}{8})\sum^{n}_{i=1}\int_{\Omega}\frac{u^2}{|x|^{4}}X^{2}_{1}(\frac{|x|}{R})%%@
X^{2}_{2 }(\frac{|x|}{R})...X^{2}_{i}(\frac{|x|}{R})dx.
\end{equation}
\end{enumerate}
Moreover, all constants in the above inequalities are best
constants.
  \end{corollary}

%We will prove the above two theorems later in this section.

\begin{theorem} \label{n-dim} Let $W_1(x)$ and $W_2(x)$ be two radial Bessel potentials on a ball %%@
$B$ of radius $R$ in $R^n$ with $n\geq 4$. Then  for %%@
all  $u \in H^{2}(B)\cap H_{0}^{1}(B)$
% \begin{eqnarray*}
%   \int_{\Omega}|\Delta u |^{2}dx&\geq&
%\frac{ C(n)(n-4)^2}{4}\int_{B}\frac{u^2}{|x|^4}dx+C(n)\beta (W_1; %%@
%R)\int_{B}W_1(x)\frac{u^2}{|x|^2}dx\\
%&&+ c(\frac{n-2a-2}{2})^2 \int_{\Omega}\frac{u^2}{|x|^{2a+2}}dx
%+c\beta (W_2; R)\int_{\Omega}W_2(x)\frac{u^2}{|x|^{2a}}dx.
% \end{eqnarray*}
%In particular, if $n\geq 5$
 \begin{eqnarray*}
   \int_{B}|\Delta u |^{2}dx &\geq& \frac{n^2(n-4)^2}{16} %%@
\int_{B}\frac{u^{2}}{|x|^{4}}dx+ \frac{n^2}{4}\beta (W_1;  R)\int_{B} %%@
W_1(x)\frac{u^{2}}{|x|^2}dx\\
&&+ \mu(\frac{n-2}{2})^2 \int_{B}\frac{u^2}{|x|^{2}}dx+\mu\beta (W_2;  %%@
R)\int_{B}W_2(x)u^2dx,
\end{eqnarray*}
% and if $n=4$, then
 % \begin{eqnarray*}
 %  \int_{B}|\Delta u |^{2}dx&\geq&
%3\beta (W_1;  R)\int_{B}W_1(x)\frac{u^2}{|x|^2}dx+ c(1-a)^2 \int_{B}\frac{u^2}{|x|^{2a+2}}dx
%+c\beta (W_2;  R)\int_{B}W_2(x)\frac{u^2}{|x|^{2a}}dx.
% \end{eqnarray*}
where $\mu$ is defined by $(\ref{mu})$.
\end{theorem}
{\bf Proof:} Here again we shall give a proof when $n\geq 5$. The case $n=4$ will be handled in %%@
the next section. We again  first use Theorem \ref{m=0} (for $n\geq 5$), then Theorem 2.15 in \cite{GM2} with the Bessel pair  $(|x|^{-2}, %%@
|x|^{-2}(\frac{(n-4)^2}{4}|x|^{-2}+W) )$, then again Theorem
\ref{main} with the Bessel pair
 $(1, (\frac{n-2}{2})^2|x|^{-2}+W)$
 to obtain
\begin{eqnarray*}
\int_{B}|\Delta u|^2 dx&\geq& \frac{n^2}{4} \int_{B}\frac{|\nabla u|^2}{|x|^2}dx+ \mu \int_{B}|\nabla u|^2 dx\\
&\geq&  \frac{n^2(n-4)^2}{16}\int_{B}\frac{u^2}{|x|^4}dx+\frac{n^2}{4} \beta (W_1;  %%@
R)\int_{B}W_1(x)\frac{u^2}{|x|^2}dx+ \mu \int_{B}|\nabla u|^2dx %%@
dx\\&\geq&
\frac{n^2(n-4)^2}{16}\int_{B}\frac{u^2}{|x|^4}dx+\frac{n^2}{4} \beta(W_1;  R) \int_{B}W_1(x)\frac{u^2}{|x|^2}dx\\
&&+ \mu(\frac{n-2}{2})^2 \int_{B}\frac{u^2}{|x|^{2}}dx +\mu
\beta(W_2;  R)\int_{B}W_2(x)u^2 dx.
 \end{eqnarray*}

\begin{theorem}\label{hrs-in} Assume $n\geq 4$ and let $W(x)$ be a Bessel potential %%@
on a ball $B$ of radius $R$ and centered at zero in $R^n$. Then the
following holds for all $u \in H^2(B) \cap H^{1}_{0}(B)$:
\begin{eqnarray}\label{hrs}
\int_{B}|\Delta u|^{2}dx &\geq&
\frac{n^{2}(n-4)^2}{16}\int_{B}\frac{u^2}{|x|^{4}}dx  \\
&&+\beta (W; R)\frac{n^2}{4} \int_{B}\frac{W(x)}{|x|^{2}}u^2 dx+
\frac{\mu^2}{R^2}||u||_{H^1_0},
\end{eqnarray}
where $\frac{\mu^2}{R^2}$ is defined by (\ref{mu}).
\end{theorem}
{\bf Proof:} Decomposing again $u \in C^{\infty}(\bar{B}_{R})$ into its spherical harmonics %%@
$\Sigma^{\infty}_{k=0}u_{k}$ where $u_{k}=f_{k}(|x|)\varphi_{k}(x)$, we calculate
\begin{eqnarray*}
\frac{1}{n\omega_{n}}\int_{R^n} |\Delta %%@
u_{k}|^2 dx&=&\int_{0}^{R}r^{n-1}(f''_{k}(r))^2dr+[n-1+2c_{k}]\int^{R}_{0}r^{%%@
n-3}(f_{k}')^2dr\\
&+&c_{k}[c_{k}+n-4]\int_{0}^{R}r^{n-5}(f_{k}(r))^2dr\\
&+&(n-1)(f_{k}'(R))^2R^{n-2m-2} \\
&\geq&\frac{n^2}{4}\int^{R}_{0}r^{n-3}(f_{k}')^2dr+ %%@
\frac{\mu^2}{R^2} \int^{R}_{0}r^{n-1}(f_{k}')^2dr\\
&+&c_{k}\int^{R}_{0}r^{n-3}(f_{k}')^2dr\\
&\geq& \frac{n^{2}(n-4)^2}{16}\int^{R}_{0}r^{n-5}(f_{k})^2dr\\
&&+\beta (W; R)\frac{n^2}{4}
\int^{R}_{0}W(r)r^{n-3}(f_{k})^2dr\\
&+&\frac{\mu^2}{R^2}
\int^{R}_{0}r^{n-1}(f_{k}')^2dr+c_{k}\frac{\mu^2}{R^2}
\int^{R}_{0}r^{n-3}(f_{k})^2
dr\\
&=&\frac{n^{2}(n-4)^2}{16n\omega_{n}}\int_{R^n}\frac{u_{k}^2}{|x|^{2m+4}}dx\\
&+&\frac{\beta (W; R)}{n\omega_{n}}(\frac{n^2}{4})
\int_{R^n}\frac{W(x)}{|x|^{2}}u_{k}^2 dx+
\frac{\mu^2}{n\omega_{n}R^2} ||u_{k}||_{W_{0}^{1,2}}.
\end{eqnarray*}
Hence (\ref{hrs}) holds. \hfill $\Box$\\

{\bf Acknowledgment:} I would like to thank Professor Nassif
Ghoussoub, my supervisor, for his valuable suggestions, constant
support, and encouragement.

 \end{document}